\documentclass[a4paper, 11pt, twoside, openany]{article}
\usepackage{amsmath,amssymb,stmaryrd, mathtools}
\usepackage{fancyhdr}
\usepackage{footmisc}
\usepackage{hyperref}
\usepackage[margin=22mm]{geometry}
\usepackage[UKenglish]{datetime}
\usepackage{graphicx}
\usepackage[usenames, dvipsnames]{color}

\numberwithin{equation}{subsection}

\newcommand{\dd}[2]{\frac{\delta #1}{\delta #2}}

\newtheorem{remark}{Remark}
\begin{document}
\begin{center}
{\Large Energy conserving upwinded compatible finite element schemes for the \\ \vspace{3mm} rotating shallow water equations}
\end{center}
\vspace{-3mm}
\hrulefill
\begin{center}
{Golo Wimmer*, Colin Cotter* and Werner Bauer*}\\
\vspace{2mm}
{\textit{* Imperial College London, **INRIA Rennes}}\\
\vspace{2mm}
\today
\end{center}
\vspace{5mm}
\begin{center}
\textbf{Abstract}
\end{center}
We present an energy conserving space discretisation of the rotating shallow water equations using compatible finite elements. It is based on an energy and enstrophy conserving Hamiltonian formulation as described in McRae and Cotter (2014), and extends it to include upwinding in the velocity and depth advection to increase stability. Upwinding for velocity in an energy conserving context was introduced for the incompressible Euler equations in Natale and Cotter (2017), while upwinding in the depth field in a Hamiltonian finite element context is newly described here. The energy conserving property is validated by coupling the spatial discretisation to an energy conserving time discretisation. Further, the discretisation is demonstrated to lead to an improved field development with respect to stability when upwinding in the depth field is included.
%
%
%
%
\section{Introduction}
The compatible finite element approach has recently been proposed as a discretisation method for numerical weather prediction \cite{COTTER20127076}. It relies on the use of so-called de Rham complexes of finite element spaces, where one space is mapped to another via differential operators. This leads to desirable stability and convergence properties and, in the context of weather prediction, further allows the use of pseudo-uniform grids on the sphere that avoid the parallel computing issues associated with the latitude-longitude grid. These issues arise due to a relatively fine mesh resolution towards the grid poles, which in turn requires more communication between mesh cells during each time step \cite{ford2013gung}. Moreover, the compatible finite element method is quite general, allowing for adaptive mesh refinement and higher-order discretisations. For this reason, the Met Office's next generation dynamical core, that is the atmosphere model's fluid dynamics component, will be based on a compatible finite element method. Further details on recent results for compatible finite elements in numerical weather prediction can be found in \cite{natale2016compatible}. \\ \\
An integral part of the governing equations used in numerical weather prediction are convection terms, and their discretisation has attracted the interest of many researchers \cite{cockburn2000development}. For finite element methods, the scheme to be used depends on the underlying finite element spaces, requiring the application of different schemes for different fields in the compatible framework. These schemes should be consistent and stable, while avoiding an excessive use of diffusion to achieve the latter. Two classical examples of such schemes are the standard upwind Discontinuous Galerkin method \cite{kuzmin2010guide} and the streamline upwind Petrov-Galerkin method \cite{brooks1982streamline} for discontinuous and continuous Galerkin finite element spaces, respectively. However, for the spaces occurring in the compatible framework in the context of numerical weather prediction, adjusted or mixed versions of the aforementioned methods may be required \cite{natale2016compatible}. \\ \\
Another aspect important for numerical weather prediction, particularly for long-term simulations, is conservation of quantities such as mass and energy. One way to ensure conservation of the latter is to discretise the governing equations within a Hamiltonian framework, where the system's Hamiltonian represents the total amount of energy. A description of the compressible Euler equations within this framework was first presented in \cite{morrison1980noncanonical} in the setting of magneto-hydrodynamics, and relies on a formulation with a Poisson bracket. This allows for the use of non-canonical (and hence in this case physical) variables \cite{shepherd1990symmetries}, and many fluid dynamical equations have been formulated in Poisson bracket form since \cite{morrison1982poisson, morrison1998hamiltonian}. Conservation of energy follows easily from this setup via the bracket's antisymmetry, and will be maintained in any space discretisation that preserves the latter. In the context of numerical weather prediction, the Poisson bracket framework has already been considered e.g. in \cite{salmon2004poisson} and \cite{gassmann2008towards} for finite difference discretisations. Further, for compatible finite element methods, it has first been considered in \cite{McRae_2014} for the rotating shallow water equations and by extension for the sphere and hemisphere with boundaries in \cite{Lee:2018aa} and \cite{BAUER2018171}, respectively. \\ \\
One way to incorporate both upwinding and energy conservation in a space discretisation is to follow the Poisson bracket framework, adding upwinding terms while ensuring that the bracket's antisymmetry is maintained. In a compatible finite element setting, this has already been achieved for potential vorticity upwinding for the rotating shallow water equations in \cite{BAUER2018171}, while upwinding for the velocity field for the incompressible Euler equations has been introduced in \cite{natale2016variational}, using a geometric approach including Lie derivatives. Further, upwinding for buoyancy has been introduced for the thermal rotating shallow water equations in \cite{ELDRED20191}. For the compressible Euler equations, energy-conserving upwinding schemes for the density and temperature fields remain to be formulated. For simplicity, we will focus on the former and revert to the rotating shallow water equations, replacing density with depth. Hence, in this paper we extend the energy-conserving formulation for the rotating shallow water equations as given in \cite{McRae_2014} to include upwinding in the depth field $D$, and further incorporate the velocity upwinding scheme of \cite{natale2016variational}. To do so, we introduce an additional operator $\mathbb{U}$ to recover the velocity field $\mathbf{u}$ from the momentum flux $\mathbf{F}$ arising in the Hamiltonian framework. The resulting discretisation is then tested for energy conservation using an energy-conserving time discretisation as introduced in \cite{cohen2011linear}, and further assessed for its qualitative field development in comparison to a version not upwinded in $D$, and a non-energy conserving version including the same type of upwinding.\\ \\
The rest of the paper is structured as follows: In section \ref{Formulation} we first review the existing, non-upwinded compatible finite element formulation of the rotating shallow water equations, and then describe the incorporation of upwinding. In section \ref{Numerical results}, we present numerical results. Finally, in section \ref{Conclusion} we review our results and discuss ongoing work.
\section{Energy conserving formulation for Rotating Shallow Water equations with upwinding terms} \label{Formulation}
In this section, we extend the energy conserving space discretisation for the rotating shallow water equations presented in \cite{McRae_2014}, by introducing an upwind formulation for the depth field, and further using the energy conserving velocity field upwinding as presented for the incompressible Euler equations in \cite{natale2016variational}. 
\subsection{Hamiltonian formulation}
To construct an energy conserving space discretisation for the rotating shallow water equations, we consider a derivation of the equations that contains a direct condition for energy conservation. It is based on the symplectic form of the Hamiltonian structure underlying fluid dynamics, and is given in terms of Poisson brackets  \cite{shepherd1990symmetries}. For Hamiltonian $H$, i.e. the system's total energy, and any functional $F$ of the dynamic variables, we have
\begin{equation}
\frac{dF}{dt} = \{F,H\}, \label{Lie-Poisson}
\end{equation}
for Poisson bracket $\{ \cdot,\cdot \}$, which is bilinear and antisymmetric. Its form depends on the choice of dynamic variables, and in the case of the rotating shallow water equations with velocity $\mathbf{u}$ and depth $D$ it is defined by
\begin{align}
\{F, H\} \coloneqq - \langle \dd{F}{\mathbf{u}}, q \dd{H}{\mathbf{u}}^\perp \rangle + \langle \nabla \cdot \dd{F}{\mathbf{u}}, \dd{H}{D} \rangle - \langle \nabla \cdot \dd{H}{\mathbf{u}}, \dd{F}{D}\rangle, \label{RSW_Poisson}
\end{align}
where $\langle\cdot,\cdot\rangle$ denotes the $L^2$ inner product over the domain $\Omega$ in consideration, $(a,b)^\perp = (b, -a)$, and $q$ denotes potential vorticity, i.e.
\begin{equation}
q = (\nabla \cdot \mathbf{u}^\perp + f)/D, \label{vorticity}
\end{equation}
for Coriolis parameter $f$. The functional derivatives are defined weakly as 
\begin{equation}
\langle \dd{F}{\mathbf{u}} , \mathbf{w} \rangle \coloneqq \lim_{\epsilon \rightarrow 0} \frac{1}{\epsilon} \big(F(\mathbf{u}+ \epsilon \mathbf{w}, D) -F(\mathbf{u}, D) \big) \hspace{2cm} \forall \mathbf{w} \in V(\Omega),
\end{equation}
for a suitable space $V(\Omega)$ to be defined, and similarly for the variation in $D$. To complete the functional equation \eqref{Lie-Poisson} for the rotating shallow water equations, we need to define the Hamiltonian $H$. It is given by
\begin{equation}
H(\mathbf{u}, D) = \frac{1}{2}  \int_\Omega (D |\mathbf{u}|^2 + g (D+b)^2)dx,\label{RSW_Hamiltonian}
\end{equation}
for gravitational acceleration $g$ and bottom profile $b$. In view of the Poisson bracket, we find that the variational derivatives are given by
\begingroup
\addtolength{\jot}{1em}
\begin{align}
&\dd{H}{\mathbf{u}} = D \mathbf{u},\\
&\dd{H}{D} = \frac{1}{2} |\mathbf{u}|^2 + g(D+b).
\end{align}
\endgroup
The usual form of the equations then follows by choosing $F = \langle \mathbf{u}, \mathbf{w} \rangle$ and $\langle D, \phi \rangle$, respectively for arbitrary test functions in $C^1(\Omega)$, noting that for the former we have $\dd{F}{D}=0$, while for the latter $\dd{F}{\mathbf{u}}=0$. Using these $F$ in \eqref{Lie-Poisson}, we recover the usual form of the rotating shallow water equations, i.e.
\begin{align}
&\mathbf{u}_t = -qD\mathbf{u} - \nabla (\frac{1}{2} |\mathbf{u}|^2 + g(D+b)) = -(\mathbf{u} \cdot \nabla) \mathbf{u} - f \mathbf{u}^\perp - g \nabla (D+b), \label{Cts_u}\\
&D_t = -\nabla \cdot (D\mathbf{u}), \label{Cts_D}
\end{align}
noting that we applied integration by parts for the second bracket term, assuming suitable boundary conditions.\\ \\
Using this framework, we find that energy conservation follows immediately due to the bracket's antisymmetry. Setting $F=H$, we find
\begin{equation}
\frac{dH}{dt} = \{H,H\} = - \{H,H\} = 0,
\end{equation}
and in particular, any space discretisation whose bracket is still antisymmetric will also satisfy conservation of energy.
\subsection{Existing discretised formulation without upwinding terms}
An energy conserving space discretisation of the rotating shallow water equations in boundary-free domains was presented in  \cite{McRae_2014}. It is based on the Hamiltonian framework as reviewed above, together with a compatible finite element discretisation. The finite element spaces for the prognostic variables $D, \mathbf{u}$ and the diagnostic potential vorticity $q$ are given by $W_2, \; W_1, \; W_0$ such that\\
\begin{equation}
W_0 \; \overset{\nabla^\perp}{\longrightarrow} \; W_1 \; \overset{\nabla \cdot}{\longrightarrow} \; W_2, \label{DeRham}
\end{equation}\\
that is the differential operators appearing in our equations map one finite element space to another. For this to be well-defined, $W_1$ is chosen to be a finite element subspace of a space for which the divergence operator is well-defined, i.e. $H_{div}(\Omega) = \{ \mathbf{u} \in L^2(\Omega; \mathbb{R}^2) \colon \nabla \cdot \mathbf{u} \in L^2(\Omega)\}$. Examples of such finite element spaces and the resulting spaces $W_0$, $W_2$ required to satisfy \eqref{DeRham} are given in \cite{McRae_2014}. For our numerical tests below, we consider the second order Brezzi-Douglas-Marini triangular finite element, which requires three point evaluations of the vector field's normal component at each of the element's edges, as well as three additional interior moments. Consequently, the potential vorticity finite element space is given by the standard triangular third order Continuous Galerkin space, while the depth space is the first order Discontinuous Galerkin space. In short, $(W_2, W_1, W_0) = (DG_1, BDM_2, CG_3)$.\\ \\
The discrete bracket is identical to the continuous one \eqref{RSW_Poisson} presented above, and complemented with an auxiliary equation for the vorticity \eqref{vorticity}, given by\\
\begin{equation}
\langle \eta, qD\rangle = -\langle \nabla^\perp \eta, \mathbf{u} \rangle + \langle \eta, f\rangle \hspace{2cm} \forall \eta \in W_0.
\end{equation}
In the discrete case, the Hamiltonian variations are given by projections into the relevant finite element spaces, i.e.
\begingroup
\addtolength{\jot}{1em}
\begin{align}
&\dd{H}{\mathbf{u}} = P_{W_1}(D \mathbf{u}) \eqqcolon \mathbf{F},\\
&\dd{H}{D} = P_{W_2}(\frac{1}{2} |\mathbf{u}|^2 + g(D+b)),
\end{align}
\endgroup
where $P_{W_1}$ denotes projection into the velocity space $W_1$, and similar for the depth space $W_2$. The resulting space-discretised weak form of the rotating shallow water equations thus reads
\begingroup
\addtolength{\jot}{1em}
\begin{align}
\begin{split}
&\langle \mathbf{w}, \mathbf{u}_t \rangle + \langle \mathbf{w}, q\mathbf{F}^\perp \rangle - \langle \nabla \cdot \mathbf{w}, \frac{1}{2} |\mathbf{u}|^2 + g(D+b)\rangle = 0 \hspace{2cm} \forall \mathbf{w} \in W_1,\\
&\langle \phi, D_t \rangle + \langle \phi , \nabla \cdot \mathbf{F} \rangle = 0 \hspace{67mm} \forall \phi \in W_2. \label{var_original_eqns}
\end{split}
\end{align}
\endgroup
Note that the divergence operator maps $\mathbf{w}$ into $W_2$, implying that the projection $P_{W_2}$ is not explicitly needed for the variation of $H$ in $D$.
\subsection{Formulation including upwinding terms}
The Poisson bracket \eqref{RSW_Poisson} leads to a natural energy-conserving space discretisation of the rotating shallow water equations in the sense that the discretised bracket is equal to the non-discretised one. However, the resulting transport schemes may produce spurious small scale features, and their stability can be improved e.g. by incorporating upwinding in $D$, $q$, or $\mathbf{u}$. In the following two subsections, we introduce a method to incorporate discontinuous Galerkin upwinding in the depth field while maintaining the bracket's antisymmetry, and show how to incorporate the energy-conserving upwinding in the velocity field as given in \cite{natale2016compatible}.
\subsubsection{DG upwinding for $D$}
In \eqref{var_original_eqns}, the depth field $D$ is advected using a discontinuous Galerkin discretisation. To improve stability, it is desirable to include an upwinding term, accounting for the total amount of depth $D$ that is flowing from one cell to another in each time step. Given an advection equation such as the depth equation
\begin{equation*}
D_t + \nabla \cdot(D \mathbf{u}) =0,
\end{equation*}
the corresponding DG upwind formulation is given by \cite{kuzmin2010guide}
\begin{equation}
\langle \phi, D_t \rangle = \langle \nabla \phi, \mathbf{u} D \rangle - \int_\Gamma [\![\phi \mathbf{u}]\!] \tilde{D} \; dS \hspace{2cm} \forall \phi \in W_2, \label{DG_upwind}
\end{equation}
where the last integral is over all mesh facets, with jump operator $[\![.]\!]$ and upwind value defined by
\begin{equation}
[\![\mathbf{v}]\!] = \mathbf{v}^+ \cdot \mathbf{n}^+ + \mathbf{v}^- \cdot \mathbf{n}^-, \hspace{2cm} \tilde{D}=
     \begin{cases}
     	\begin{split}
       	D^+ \; \; \; \; \text{if } \mathbf{u} \cdot \mathbf{n}^+ < 0\\
	D^- \; \; \; \; \text{otherwise}, \label{def_upwind_D} \; \; \;
	\end{split}
     \end{cases}
\end{equation}
noting that the two sides of each mesh facet are arbitrarily denoted by + and - (and hence $\mathbf{n}^+ = -\mathbf{n}^-$). If the solution of a given problem is smooth, we find that the facet term in \eqref{DG_upwind} vanishes as we refine our discretisation's resolution, showing that the incorporation of upwinding is still consistent.\\ \\
In the context of the energy-conserving discretisation reviewed above, note that in \eqref{DG_upwind} the velocity and depth fields $\mathbf{u}$, $D$ are separated in the facet integral, while they appear implicitly in the flux projection $\mathbf{F}= P_{W_1}(D \mathbf{u})$ in the discretised equation set  \eqref{var_original_eqns}. Possible alternative forms of the facet integral including $\mathbf{F}$, such as
\begin{align}
&\int_\Gamma [\![\phi \mathbf{F}/D]\!] \tilde{D} \; dS, \label{no mass conservation} \\
&\int_\Gamma [\![\phi \mathbf{F}]\!] \tilde{D} \; dS, \label{not consistent}\\ 
&\int_\Gamma [\![\phi]\!]  \text{avg}(\mathbf{F}/D)\tilde{D} \; dS, \label{bad fields avg}
\end{align}
are not mass-conserving (\eqref{no mass conservation}), or can be shown to produce no improvement with respect to stability (\eqref{not consistent}) or unstable fields (\eqref{bad fields avg}). To incorporate the standard upwinding \eqref{DG_upwind} in an energy conserving discretisation, we will need to introduce an additional operator to avoid the occurrence of projections due to the discrete Hamiltonian variations. In particular, we introduce an operator to recover our velocity $\mathbf{u}$ from $\mathbf{F}$, arising due to the Hamiltonian variation with respect to $\mathbf{u}$. Hence, define $\mathbb{U}$ implicitly by
\begin{equation}
\hspace{2cm} \mathbb{U}(D, \mathbf{F}) \colon  W_2 \times W_1 \longrightarrow W_1 \; \; \; \; \text{such that} \; \; \; \; \langle D \mathbf{v}, \mathbb{U} \rangle = \langle \mathbf{v}, \mathbf{F} \rangle \hspace{15mm} \forall \mathbf{v} \in W_1. \label{bbU_eqn}
\end{equation}
Note that in the continuous sense, $\mathbb{U}(D, \mathbf{F})$ corresponds to division of $\mathbf{F}$ by $D$. Since in our shallow water setting $D > 0$, we find that $\mathbb{U}$ is well-defined and further that if $\mathbf{F}$ is any function in $W_1$ of flux form $P_{W_1}(D \mathbf{u})$, then
\begin{equation}
\mathbb{U}(D, \mathbf{F}) = \mathbf{u} \label{bbU_pointwise}
\end{equation}
pointwise.\\ \\
Using $\mathbb{U}$, we are in a position to alter bracket \eqref{RSW_Poisson} introducing upwinding terms for $D$, and arrive at a Poisson bracket of form
\begin{align}
\{F, H\} = -\langle \frac{\delta F}{\delta \mathbf{u}}, q \frac{\delta H}{\delta \mathbf{u}}^\perp \rangle & - \langle D \; \mathbb{U}\big(D, \frac{\delta F}{\delta \mathbf{u}}\big), \nabla \frac{\delta H}{\delta D}\rangle  + \int_\Gamma [\![\frac{\delta H}{\delta D} \mathbb{U}\big(D, \frac{\delta F}{\delta \mathbf{u}}\big)]\!] \tilde{D} \; dS\\ \nonumber \\
& +  \langle D \; \mathbb{U}\big(D, \frac{\delta H}{\delta \mathbf{u}}\big), \nabla \frac{\delta F}{\delta D}\rangle - \int_\Gamma [\![\frac{\delta F}{\delta D} \mathbb{U}\big(D, \frac{\delta H}{\delta \mathbf{u}}\big)]\!] \tilde{D} \; dS, \label{var_Gv1_bracket}
\end{align}
noting that we also had to introduce upwinding for the corresponding antisymmetric term of the momentum equation (second term in \ref{RSW_Poisson}) to maintain antisymmetry. Further, note that \eqref{var_Gv1_bracket} corresponds to the standard upwinding formulation \eqref{DG_upwind}, provided that \eqref{bbU_pointwise} holds.
\subsubsection{Upwinding in $\mathbf{u}$}
It is possible to also apply upwinding in $q$, as done in \cite{BAUER2018171}. Alternatively, we can aim to improve stability for the evolution of the velocity field $\mathbf{u}$ by considering the momentum equation purely in $\mathbf{u}$ and $D$, recalling that
\begin{equation}
q = (\nabla^\perp \cdot \mathbf{u} + f)/D.
\end{equation}
Using this form of $q$, the term in the weak equations (cf \eqref{var_original_eqns}) corresponding to the vorticity part of the bracket reads
\begin{equation}
\langle \mathbf{w}, q\mathbf{F}^\perp \rangle = \langle \mathbf{w}, qD\mathbf{u}^\perp \rangle = \langle \mathbf{w},( \nabla^\perp \cdot \mathbf{u})\mathbf{u}^\perp \rangle + \langle \mathbf{w}, f\mathbf{u}^\perp \rangle, \label{no_q}
\end{equation}\\
noting that for the purpose of presenting the velocity upwinding scheme, we temporarily ignore the projection in $\mathbf{F}$.  We can then replace this form by an upwind formulation for $\mathbf{u}$, e.g. as introduced in \cite{natale2016compatible}:
\begin{equation}
- \langle \nabla^\perp (\mathbf{w} \cdot \mathbf{u}^\perp), \mathbf{u} \rangle + \int_\Gamma [\![\mathbf{w} \cdot \mathbf{u}^\perp]\!] \; \mathbf{n}^\perp \cdot \tilde{\mathbf{u}} \; dS + \langle \mathbf{w}, f\mathbf{u}^\perp \rangle. \label{nonvar_u_upwind}
\end{equation}
To formulate a Poisson bracket that leads to a set of governing equations with this form of upwinding, we need to use our velocity recovering operator $\mathbb{U}$ again to avoid projections of form $P_{W_1}(D\mathbf{u})$. Note that in terms of functionals $F$ and $H$, the formulation of \eqref{no_q} without $q$ and its upwinded extension \eqref{nonvar_u_upwind} read
\begin{align}
\langle \dd{F}{\mathbf{u}}, q\dd{H}{\mathbf{u}}^\perp \rangle = &\langle \dd{F}{\mathbf{u}}, ( \nabla^\perp \cdot \mathbf{u})\big(\dd{H}{\mathbf{u}}/D\big)^\perp \rangle + \langle \dd{F}{\mathbf{u}}, f \big(\dd{H}{\mathbf{u}}/D\big)^\perp \rangle\\
\rightarrow & - \langle \nabla^\perp (\dd{F}{\mathbf{u}} \cdot \big(\dd{H}{\mathbf{u}}/D\big)^\perp), \mathbf{u} \rangle + \int_\Gamma [\![\dd{F}{\mathbf{u}} \cdot \big(\dd{H}{\mathbf{u}}/D\big)^\perp]\!] \; \mathbf{n}^\perp \cdot \tilde{\mathbf{u}} \; dS \nonumber \\ 
&+ \langle\dd{F}{\mathbf{u}}, f\big(\dd{H}{\mathbf{u}}/D\big)^\perp \rangle.
\end{align}
We find that this time, since each of the advection bracket terms are by themselves antisymmetric, we have to apply $\mathbb{U}$ twice in each integral: once for $\dd{H}{\mathbf{u}}$ to recover the upwinding velocity $\mathbf{u}$, and once for $\dd{F}{\mathbf{u}}$ to maintain antisymmetry. Note that as $\mathbb{U}$ corresponds to division by $D$, we introduce, for the purpose of consistency, an additional depth term $D$ wherever $\mathbb{U}$ is not applied to a velocity space element of flux form $P_{W_1}(D\mathbf{u})$. Thus, the above upwinding form \eqref{nonvar_u_upwind} is given in the Hamiltonian variational setting by
\begingroup
\addtolength{\jot}{1em}
\begin{align}
\langle \nabla^\perp \Big(D \mathbb{U}\big(D, \dd{F}{\mathbf{u}}\big) \cdot \mathbb{U}\big(D, \dd{H}{\mathbf{u}}\big)^\perp \Big), \mathbf{u} \rangle &- \int_\Gamma [\![D \mathbb{U}\big(D, \dd{F}{\mathbf{u}}\big) \cdot \mathbb{U}\big(D, \dd{H}{\mathbf{u}}\big)^\perp]\!] \mathbf{n}^\perp \cdot \tilde{\mathbf{u}} \; dS\\
&- \langle D \mathbb{U}\big(D, \dd{F}{\mathbf{u}}\big), f \mathbb{U}\big(D, \dd{H}{\mathbf{u}}\big)^\perp\rangle.
\end{align}
\endgroup
Altogether, the full Poisson bracket including upwinding for depth and velocity advection is hence given by
\begingroup
\addtolength{\jot}{1em}
\begin{align}
\{F, H\} &\coloneqq  \langle \nabla^\perp \big(D \mathbb{U}\big(D, \dd{F}{\mathbf{u}}\big) \cdot \mathbb{U}\big(D, \dd{H}{\mathbf{u}}\big)^\perp \big), \mathbf{u} \rangle - \int_\Gamma [\![D \mathbb{U}\big(D, \dd{F}{\mathbf{u}}\big) \cdot \mathbb{U}\big(D, \dd{H}{\mathbf{u}}\big)^\perp]\!] \mathbf{n}^\perp \cdot \tilde{\mathbf{u}} \; dS \label{var_Gv2_bracket_u_a}\\
&- \langle D \; \mathbb{U}\big(D, \frac{\delta F}{\delta \mathbf{u}}\big), \nabla \frac{\delta H}{\delta D}\rangle + \int_\Gamma [\![\frac{\delta H}{\delta D} \mathbb{U}\big(D, \frac{\delta F}{\delta \mathbf{u}}\big)]\!] \tilde{D} \; dS - \langle D \mathbb{U}\big(D, \dd{F}{\mathbf{u}}\big), f \mathbb{U}\big(D, \dd{H}{\mathbf{u}}\big)^\perp\rangle \label{var_Gv2_bracket_u_f}\\
&+  \langle D \; \mathbb{U}\big(D, \frac{\delta H}{\delta \mathbf{u}}\big), \nabla \frac{\delta F}{\delta D}\rangle - \int_\Gamma [\![\frac{\delta F}{\delta D} \mathbb{U}\big(D, \frac{\delta H}{\delta \mathbf{u}}\big)]\!] \tilde{D} \; dS. \label{var_Gv2_bracket_D_a}
\end{align}
\endgroup
Checking for antisymmetry, we find that the terms in line \eqref{var_Gv2_bracket_u_a} as well as the Coriolis term are antisymmetric by themselves (due to the perpendicular $(a, b)^\perp = (b, -a)$), while the first and second terms respectively in \eqref{var_Gv2_bracket_u_f} and \eqref{var_Gv2_bracket_D_a} form antisymmetric pairs.\\ \\
In view of the time discretisations to follow and for ease of notation, we rewrite the above bracket as
\begin{equation}
\{F, H\}  = L_{\mathbb{U}(D, \frac{\delta H}{\delta \mathbf{u}})}\big( \mathbf{u}; D \mathbb{U}(D, \dd{F}{\mathbf{u}})\big) + F_{(D, \mathbb{U}(D, \dd{H}{\mathbf{u}}), \dd{H}{D})}\big(\mathbb{U}(D, \dd{F}{\mathbf{u}})\big) + L^D_{\mathbb{U}(D, \frac{\delta H}{\delta \mathbf{u}})}\big(D; \dd{F}{D}\big), \label{var_Gv2_bracket}
\end{equation}
with velocity advection operator $L$ corresponding to \eqref{var_Gv2_bracket_u_a}, forcing operator $F$ to \eqref{var_Gv2_bracket_u_f}, and depth advection operator $L^D$ to \eqref{var_Gv2_bracket_D_a}. The choice of notation for advection is based on the usual notation for the Lie derivative $\pounds_u$, which in the case of advecting velocity field $\mathbf{v}$ and advected velocity field $\mathbf{u}$ is given by
\begin{equation}
\pounds_{\mathbf{v}}(\mathbf{u}) = (\nabla^\perp \cdot \mathbf{u})\mathbf{v}^\perp + \nabla (\mathbf{v} \cdot \mathbf{u}).
\end{equation}
More specifically, comparing this to \eqref{nonvar_u_upwind}, we see that $L$ corresponds to the divergence-free part of velocity advection (in the Lie derivative sense, noting that in our case the divergence part is contained in $\dd{H}{D}$). Finally, the choice of notation for forcing $F$ is simply to resemble $L$, in that the lowered terms indicate acting fields. Note that since $F$ does not explicitly contain the velocity field it is acting on, its only explicit argument in this notation is the test function $\mathbb{U}(D, \dd{F}{\mathbf{u}})$.
\begin{remark}
Comparing the advection form \eqref{nonvar_u_upwind}, which was derived in a non-energy conserving context, with the Poisson bracket based form introduced above, we find
\begin{equation}
L_{\mathbf{u}} \big(\mathbf{u}; \mathbf{w}\big) \; \; \; \longrightarrow \; \; \; L_{\mathbb{U}(D, \mathbf{F})}( \mathbf{u}; D \mathbb{U}(D, \mathbf{w})\big),
\end{equation}
that is we replaced the advecting velocity $\mathbf{u}$ by the flux-recovered velocity $\mathbb{U}(D, \mathbf{F})$, while the advected velocity remains unchanged.  Further, the test function $\mathbf{w}$ is replaced by a discrete multiplication and division of $\mathbf{w}$ by $D$. Similarly, for advection in $D$ the advection velocity is now also given by $\mathbb{U}(D, \mathbf{F})$ compared to standard DG upwinding \eqref{DG_upwind}:
\begin{equation}
L^D_{\mathbf{u}} (D; \phi) \; \; \; \longrightarrow \; \; \; L^D_{\mathbb{U}(D, \mathbf{F})}(D; \phi).
\end{equation}
\end{remark}
\section{Numerical results} \label{Numerical results}
In the previous section we introduced bracket \eqref{var_Gv2_bracket_u_a} - \eqref{var_Gv2_bracket_D_a}, which is based on the variational scheme \eqref{var_original_eqns} as given in \cite{McRae_2014} and extends it to include upwinding in the depth and velocity fields. To demonstrate conservation of energy, we additionally use an energy conserving time discretisation, thus expecting energy conservation to machine precision. Before moving on to the test cases, we review the time discretisation as well as the solver scheme for the resulting nonlinear system of equations.
\subsection{Energy conserving time discretisation}
The time integrator is given by a Poisson integrator, which was introduced in \cite{cohen2011linear} and ensures conservation of higher degree Hamiltonians. This includes the Hamiltonian corresponding to the total energy of the shallow water equations, which is cubic. To use the integrator, we follow \cite{BAUER2018171} and exploit the fact that our Poisson bracket formulation \eqref{Lie-Poisson} can be written as a system of form
\begin{equation}
\dot{\mathbf{z}} = J(\mathbf{z}) \dd{}{\mathbf{z}} H(\mathbf{z}), \label{J_system}
\end{equation}
for unknown $\mathbf{z} = (\mathbf{u}, D)$, Hamiltonian $H$ and a skew-symmetric transformation $J$ determined by the Poisson bracket via the relation
\begin{equation}
\{F, H\} = \langle \dd{F}{\mathbf{z}},  J(\mathbf{z}) \dd{H}{\mathbf{z}} \rangle.
\end{equation}
A Hamiltonian conserving time integrator for a system of ODEs of form \eqref{J_system} is then given by
\begin{equation}
\mathbf{z}^{n+1} = \mathbf{z}^n + \Delta t J\Big(\frac{\mathbf{z}^{n+1} + \mathbf{z}^n}{2}\Big) \Big(\overline{\dd{H}{\mathbf{u}}}, \overline{\dd{H}{D}}\Big), \label{J_system_discretised}
\end{equation}
with time-averaged Hamiltonian given by
\begin{equation}
\overline{\dd{H}{\mathbf{u}}} \coloneqq \int_0^1 \dd{}{\mathbf{u}} H(\mathbf{z}^n + s(\mathbf{z}^{n+1} - \mathbf{z}^n)) ds,
\end{equation}
and similar for $\overline{\dd{H}{D}}$. In our case, we can integrate the time-averaged Hamiltonians and find
\begin{align}
\overline{\dd{H}{\mathbf{u}}} =& \frac{1}{3} P_{W_1} \big(D^n \mathbf{u}^n + \frac{1}{2} D^n \mathbf{u}^{n+1} + \frac{1}{2} D^{n+1} \mathbf{u}^n + D^{n+1} \mathbf{u}^{n+1} \big),\\
\overline{\dd{H}{D}} =& P_{W_2} \big(\frac{1}{6} (|\mathbf{u}^n|^2 + \mathbf{u}^n \cdot \mathbf{u}^{n+1} + |\mathbf{u}^{n+1}|^2) +  g (\frac{1}{2}(D^n + D^{n+1}) + b) \big).
\end{align}
\begin{remark}
Since $\overline{\dd{H}{\mathbf{u}}}$ is not of a simple flux form $D \mathbf{u}$ anymore, we find that the pointwise relation \eqref{bbU_pointwise} for our velocity recovering operator does not hold anymore  for this time scheme, and we have to revert to the defining relation of $\mathbb{U}$, now given by
\begin{equation}
\langle \frac{1}{2}(D^n + D^{n+1}) \mathbf{v}, \mathbb{U} \rangle = \langle \mathbf{v}, \mathbf{F} \rangle \hspace{3cm} \forall \mathbf{v} \in W_1. \label{poisson_bbU}
\end{equation}
Note that in accordance with \eqref{J_system_discretised}, we choose a midpoint time average for $D$ here since in view of the Poisson system \eqref{J_system}, this relation is part of the transformation $J$. For $\overline{\dd{H}{\mathbf{u}}}$, the above relation \eqref{poisson_bbU} is hence given by\\
\begin{equation}
\langle \frac{1}{2}(D^n + D^{n+1}) \mathbf{v}, \mathbb{U} \rangle = \frac{1}{3} \langle \mathbf{v}, D^n \mathbf{u}^n + \frac{1}{2} D^n \mathbf{u}^{n+1} + \frac{1}{2} D^{n+1} \mathbf{u}^n + D^{n+1} \mathbf{u}^{n+1} \rangle \hspace{1cm} \forall \mathbf{v} \in W_1, \label{poisson_bbU_H}
\end{equation}\\
replacing the pointwise version \eqref{bbU_pointwise}.
\end{remark}
Writing $\bar{D} = \frac{1}{2}(D^n + D^{n+1})$, $\bar{\mathbf{u}} = \frac{1}{2}(\mathbf{u}^n + \mathbf{u}^{n+1})$, and $\bar{\mathbb{U}}$ for the solution to \eqref{poisson_bbU_H}, we arrive at a fully discretised set of nonlinear equations of form
\begin{align}
\langle \mathbf{w}, \mathbf{u}^{n+1}-\mathbf{u}^n \rangle =& \Delta t \Big(L_{\bar{\mathbb{U}}}\big( \bar{\mathbf{u}}; \bar{D} \mathbb{U}(\bar{D}, \mathbf{w})\big) + F_{(\bar{D}, \bar{\mathbb{U}}, \overline{\dd{H}{D}} )}\big(\mathbb{U}(\bar{D}, \mathbf{w})\big) \Big) &\forall \mathbf{w} \in W_1, \label{discrete_u_eqn}\\
\langle \phi, D^{n+1} - D^n \rangle =&  \Delta t L^D_{\bar{\mathbb{U}}}(\bar{D}; \phi) &\forall \phi \in W_2, \label{discrete_D_eqn}
\end{align}
to be solved for $\mathbf{u}^{n+1}$, $D^{n+1}$. Note that the Poisson integrator also requires upwinding using $\bar{\mathbf{u}}$, as the upwinded $\tilde{D}$ is part of the transformation $J$ (cf \eqref{def_upwind_D}).
\subsubsection{Nonlinear solver}
In this subsection, we briefly describe the scheme used for finding a solution for the nonlinear system of equations \eqref{discrete_u_eqn} - \eqref{discrete_D_eqn}. We revert to a Picard iteration scheme, starting from the governing equations in a residual formulation:
\begin{equation}
0 = \mathbf{R}(\mathbf{z}^{n+1}; (\mathbf{w}, \phi)^T) = (R_\mathbf{u}, R_D)^T,
\end{equation}
for $R_\mathbf{u}, \; R_D$ defined as the difference of the left-hand and right-hand sides of  \eqref{discrete_u_eqn} and \eqref{discrete_D_eqn} respectively. In the iteration scheme, we aim to find the next time value $\mathbf{z}^{n+1,k+1}$ given the old value $\mathbf{z}^n$ and the latest guess for the next time value $\mathbf{z}^{n+1,k}$. Using an increment $\delta \mathbf{z} \coloneqq \mathbf{z}^{n+1,k+1} - \mathbf{z}^{n+1,k}$, we find
\begin{align}
\begin{split}
0 = \mathbf{R}(\mathbf{z}^{n+1}; (\mathbf{w}, \phi)^T) &\approx \mathbf{R}(\mathbf{z}^{n+1, k+1}; (\mathbf{w}, \phi)^T) \\
&=  \mathbf{R}(\mathbf{z}^{n+1,k} +\delta \mathbf{z};(\mathbf{w}, \phi)^T) \\
&= \mathbf{R}(\mathbf{z}^{n+1,k}; (\mathbf{w}, \phi)^T)  + \langle \frac{\delta \mathbf{R}}{\delta \mathbf{z}} \; \delta \mathbf{z}, (\mathbf{w}, \phi)^T\rangle + O(\| \delta \mathbf{z} \|^2)\\
&\approx \mathbf{R}(\mathbf{z}^{n+1,k};(\mathbf{w}, \phi)^T) + \langle \frac{\delta \mathbf{R}}{\delta \mathbf{z}} \; \delta \mathbf{z}, (\mathbf{w}, \phi)^T \rangle,
\end{split}
\end{align}
and hence
\begin{equation}
- \mathbf{R}(\mathbf{z}^{n+1,k};(\mathbf{w}, \phi)^T) \approx \langle \frac{\delta \mathbf{R}}{\delta \mathbf{z}} \; \delta \mathbf{z}, (\mathbf{w}. \phi)^T\rangle. \label{Newton}
\end{equation}
To treat the right-hand side variational derivative, we first simplify by considering the residual $\mathbf{R}'$ derived from a weak form of the continuous equations \eqref{Cts_u} - \eqref{Cts_D} instead of $\mathbf{R}$, thus avoiding projections introduced in the energy-conserving framework. Further, we revert to a Picard iteration scheme by linearising over a background state given by $(\mathbf{u}, D)=(\mathbf{0}, h)$ for reference height $h$. Altogether, we then arrive at a right-hand side of form
 \begin{align}
\langle \frac{\delta \mathbf{R}}{\delta \mathbf{z}} \; \delta \mathbf{z}, (\mathbf{w}, \phi)^T\rangle\approx  
\begin{pmatrix}
\langle \delta \mathbf{u}, \mathbf{w} \rangle +\frac{\Delta t}{2} \langle f \delta \mathbf{u}^\perp, \mathbf{w} \rangle- \frac{\Delta t}{2} \langle gD, \nabla \cdot \mathbf{w} \rangle\\
\langle \delta D, \phi \rangle + \frac{\Delta t}{2} \langle h \nabla \cdot \mathbf{\delta u}, \phi \rangle \label{Linearised_Jacobian}
\end{pmatrix}.
\end{align}
\begin{remark} \label{bbU_test}
In order to consider the left-hand side of \eqref{Newton}, we need to find a way to treat the velocity recovery operator applied to test functions, i.e. $\mathbb{U}(\bar{D}, \mathbf{w})$, in the fully discretised momentum equation \eqref{discrete_u_eqn}. Noting that $\mathbb{U}$ corresponds to a discrete division by $\bar{D}$, this can be done by using test functions weighted by $\bar{D}$. That is, to solve for general $G\big(\mathbb{U}(\bar{D}, \mathbf{w})\big)$, we can find $\mathbf{u}$ such that
\begin{equation}
\langle \mathbf{u}, \bar{D} \mathbf{v} \rangle = G(\mathbf{v}) \hspace{2cm} \forall \mathbf{v} \in W_1. \label{remove_bbU}
\end{equation}
Then in particular, for any given test function $\mathbf{w}$, we have
\begin{equation}
G\big(\mathbb{U}(\bar{D}, \mathbf{w})\big) = \langle \mathbf{u}, \bar{D} \mathbb{U}(\bar{D}, \mathbf{w}) \rangle = \langle \mathbf{u}, \mathbf{w} \rangle,
\end{equation}
where we used \eqref{remove_bbU} for the first equality and the time-discrete defining relation \eqref{poisson_bbU} of $\mathbb{U}$ for the second one, noting that here $\mathbf{w}$ plays the role of $\mathbf{F}$ in the defining relation, while $\mathbf{u}$ corresponds to a particular choice of test function $\mathbf{v}$.
\end{remark}
Finally, the residual $-\mathbf{R}(\mathbf{z}^{n+1,k};(\mathbf{w}, \phi)^T)$ can be calculated directly using the discretised equations \eqref{discrete_u_eqn} and \eqref{discrete_D_eqn}. For the momentum equation, we first find forcing and advection velocities $\mathbf{u}^f, \mathbf{u}^a$ given by
\begin{align}
\langle \mathbf{u}^f, \bar{D} \mathbf{v} \rangle = & \Delta t F_{(\bar{D}, \bar{\mathbb{U}}, \overline{\dd{H}{D}} )}(\mathbf{v}) &\forall \mathbf{v} \in W_1, \label{poisson_forcing}\\
\langle \mathbf{u}^a, \bar{D} \mathbf{v} \rangle = & \Delta t L_{\bar{\mathbb{U}}}\big( \bar{\mathbf{u}}; \bar{D} \mathbf{v}\big) +  \langle \mathbf{u}^n, \bar{D} \mathbf{v} \rangle &\forall \mathbf{v} \in W_1, \label{poisson_advection}
\end{align}
that is $\mathbf{u}^f$ corresponds to the additional velocity induced by forcing, while $\mathbf{u}^a$ corresponds to $\mathbf{u}^n$ after advection. Next, given that \eqref{poisson_advection} and \eqref{poisson_forcing} hold, we find that in particular they hold for $\mathbf{v} = \mathbb{U}(\bar{D}, \mathbf{w}) \in W_1$ (as described in remark \ref{bbU_test}), so that for any $\mathbf{w} \in W_1$, \eqref{discrete_u_eqn} can be reformulated to\\
\begin{align}
R_{\mathbf{u}}(\mathbf{z}^{n+1,k}; \mathbf{w}) =&  \langle \mathbf{w}, \mathbf{u}^{n+1,k}-\mathbf{u}^n \rangle - \Delta t \Big(L_{\bar{\mathbb{U}}}\big( \bar{\mathbf{u}}; \bar{D} \mathbb{U}(\bar{D}, \mathbf{w})\big) + F_{(\bar{D}, \bar{\mathbb{U}}, \overline{\dd{H}{D}} )}\big(\mathbb{U}(\bar{D}, \mathbf{w})\big) \Big)\\
=& \langle \mathbf{w}, \mathbf{u}^{n+1,k} \rangle - \langle \bar{D} \mathbb{U}(\bar{D}, \mathbf{w}), \mathbf{u}^{n} \rangle \nonumber\\ 
&- \Delta t L_{\bar{\mathbb{U}}}\big( \bar{\mathbf{u}}; \bar{D} \mathbb{U}(\bar{D}, \mathbf{w})\big) - \Delta t F_{(\bar{D}, \bar{\mathbb{U}}, \overline{\dd{H}{D}} )}\big(\mathbb{U}(\bar{D}, \mathbf{w})\big)\\
=& \langle \mathbf{w}, \mathbf{u}^{n+1,k} \rangle - \langle \bar{D} \mathbb{U}(\bar{D}, \mathbf{w}), \mathbf{u}^a \rangle - \langle \bar{D} \mathbb{U}(\bar{D}, \mathbf{w}), \mathbf{u}^f \rangle\\
=& \langle \mathbf{w}, \mathbf{u}^{n+1,k} - \mathbf{u}^a - \mathbf{u}^f \rangle,
\end{align}
noting that we used the time-discrete defining relation \eqref{poisson_bbU} of $\mathbb{U}$ for the flux mass terms (i.e. $\langle \bar{D} \mathbb{U}(\bar{D}, \mathbf{w}), \cdot \rangle$) of $\mathbf{u}^n$, $\mathbf{u}^a$ and $\mathbf{u}^f$. Similarly, for the residual in $D$ we find
\begin{equation}
R_D(\mathbf{z}^{n+1,k}; \phi) = \langle{\phi, D^{n+1,k}  - D^a}\rangle,
\end{equation}
where $D^a$ corresponds to $D^n$ after advection and is solved for analogously to $\mathbf{u}^a$ in \eqref{poisson_advection}. Note that $\mathbf{u}^a + \mathbf{u}^f$ and $D^a$ can be seen as a guess for the next iteration value $k+1$. Further, note that the residual $\mathbf{R}$ on the left-hand side of \eqref{Newton} is explicit in that it depends on the known values $\mathbf{z}^n$ and $\mathbf{z}^{n+1,k}$ only. To increase the scheme's robustness, we can instead also solve for a more implicit residual system of form
\begin{align}
\langle \mathbf{u}^a, \bar{D} \mathbf{v} \rangle = & \Delta t L_{\bar{\mathbb{U}}}\big( \hat{\mathbf{u}}; \bar{D} \mathbf{v}\big) +  \langle \mathbf{u}^n, \bar{D} \mathbf{v} \rangle &\forall \mathbf{v} \in W_1,\\
\langle D^a, \phi \rangle = & \Delta t L^D_{\bar{\mathbb{U}}} \big(\hat{D}, \phi\big) + \langle D^n, \phi \rangle & \forall \phi \in W_2,
\end{align}
where we replaced the known advected time-averages $\bar{\mathbf{u}}, \bar{D}$ by implicit averages
\begin{equation}
\hat{\mathbf{u}} \coloneqq (\mathbf{u}^n + \mathbf{u}^a)/2, \; \; \; \hat{D} \coloneqq (D^n + D^a)/2.  
\end{equation}
Further, $D^a$ now also appears in the forcing term via a modified variation $\overline{\dd{H}{D}}$, given by
\begin{equation}
\overline{\dd{H}{D}} = P_{W_2} \big(\frac{1}{6} (|\mathbf{u}^n|^2 + \mathbf{u}^n \cdot \mathbf{u}^{n+1,k} + |\mathbf{u}^{n+1,k}|^2) +  g (\frac{1}{2}(D^n + D^a) + b) \big),
\end{equation}
i.e. we replaced the known $D^{n+1, k}$ by $D^a$. Note that this implicit setup constitutes a different Picard iteration scheme, which, however, can be shown to converge to the same solution as the more explicit version \eqref{poisson_forcing} - \eqref{poisson_advection}. For the test cases used below, we will use this form for the higher resolution Galewsky test case.\\ \\
This completes the fully energy-conserving scheme, with numerical test results presented in subsection \ref{Test cases}. The calculations are performed using the automated finite element toolkit Firedrake\footnote{see \url{http://firedrakeproject.org}} \cite{rathgeber2016firedrake}, using a hybridised solver to solve for the updates $\delta \mathbf{z}$.
\begin{remark} \label{non_conservative_time_disc}
A simpler time discretisation, albeit non-energy conserving, would be to use a midpoint rule for both $\mathbf{u}$ and $D$. If we further time-discretise the Hamiltonian variation in $\mathbf{u}$ as $\dd{H}{\mathbf{u}} = P_{W_1} (\bar{D} \bar{\mathbf{u}})$ and the velocity recovery operator as before (\eqref{poisson_bbU}), we find that the pointwise relation \eqref{bbU_pointwise} holds again. The resulting left-hand side equations to be solved then read
\begin{align}
\langle \mathbf{u}^f, \bar{D} \mathbf{v} \rangle = & \Delta t F_{(\bar{D}, \bar{\mathbf{u}}, \overline{\dd{H}{D}} )}(\mathbf{v}) &\forall \mathbf{v} \in W_1,\\
\langle \mathbf{u}^a, \bar{D} \mathbf{v} \rangle = & \Delta t L_{\bar{\mathbf{u}}}\big( \bar{\mathbf{u}}; \bar{D} \mathbf{v}\big) +  \langle u^n, \bar{D} \mathbf{v} \rangle &\forall \mathbf{v} \in W_1,\\
\langle D^a , \phi \rangle = &\Delta t L_{\bar{\mathbf{u}}}^D(\bar{D}, \phi) &\forall \phi \in W_2.
\end{align} 
In particular, we find that in this case only one projection, i.e. $\overline{\dd{H}{D}}$, occurs. Since it is a projection into the discontinuous Galerkin space $W_2$, the additional cost of calculating it in each Picard iteration is low. In contrast, for no upwinding in $D$ (as described in \eqref{u-ad_flux} below), this time discretisation leads to an additional Hamiltonian variation in $\mathbf{u}$ rather than $D$ (appearing in the depth advection term), leading to a higher increase in computational cost as the underlying velocity space $W_1$ is not a DG space (i.e. contains nodes on cell boundaries, leading to a non-block diagonal matrix).
\end{remark}
\subsection{Comparison to other discretisations}
To test the newly introduced upwinding in $D$ for the energy conservation as well as the qualitative field development, we compare our upwinded formulation \eqref{var_Gv2_bracket}, with one that includes upwinding in $\mathbf{u}$ only, and one that does not conserve energy. The former is given by a Poisson bracket with velocity terms equal to \eqref{var_Gv2_bracket_u_a}, but velocity forcing and depth advection terms of non-upwinded form, i.e.
\begin{equation}
\{F, H\} = L_{\mathbb{U}(D, \frac{\delta H}{\delta \mathbf{u}})}\big( \mathbf{u}; D \mathbb{U}(D, \dd{F}{\mathbf{u}})\big) + \langle \nabla \cdot \dd{F}{\mathbf{u}}, \dd{H}{D} \rangle - \langle \nabla \cdot \dd{H}{\mathbf{u}}, \dd{F}{D}\rangle. \label{u-ad_flux}
\end{equation}
As a non-energy conserving discretisation, we use standard DG upwinding for the depth field $D$, velocity upwinding of form \eqref{nonvar_u_upwind}, and forcing equal to
\begin{equation}
\langle \nabla \cdot \mathbf{w},  \frac{1}{2} |\mathbf{u}|^2 + g(D+b) \rangle.
\end{equation}
Note that in order to apply the same time discretisation to the non-energy conserving spatial discretisation, we need to rewrite the latter in terms of Hamiltonian variations and a bracket. It is given by
\begin{equation}
\{F, H\} = L_{\mathbb{U}(D, \frac{\delta H}{\delta \mathbf{u}})}\big( \mathbf{u}; \dd{F}{\mathbf{u}}\big) + L^D_{\mathbb{U}(D, \frac{\delta H}{\delta \mathbf{u}})}\big(D; \dd{F}{D}\big)+ \langle \nabla \cdot \dd{F}{\mathbf{u}}, \dd{H}{D} \rangle, \label{Non_Var}
\end{equation}
and we find that as opposed to the energy-conserving upwinded version, the velocity advection operator is not antisymmetric in itself, and the depth advection and velocity forcing are not antisymmetric to each other.
\subsection{Test cases} \label{Test cases}
\begin{figure}[ht]
\begin{center}
\includegraphics[width=0.49\textwidth]{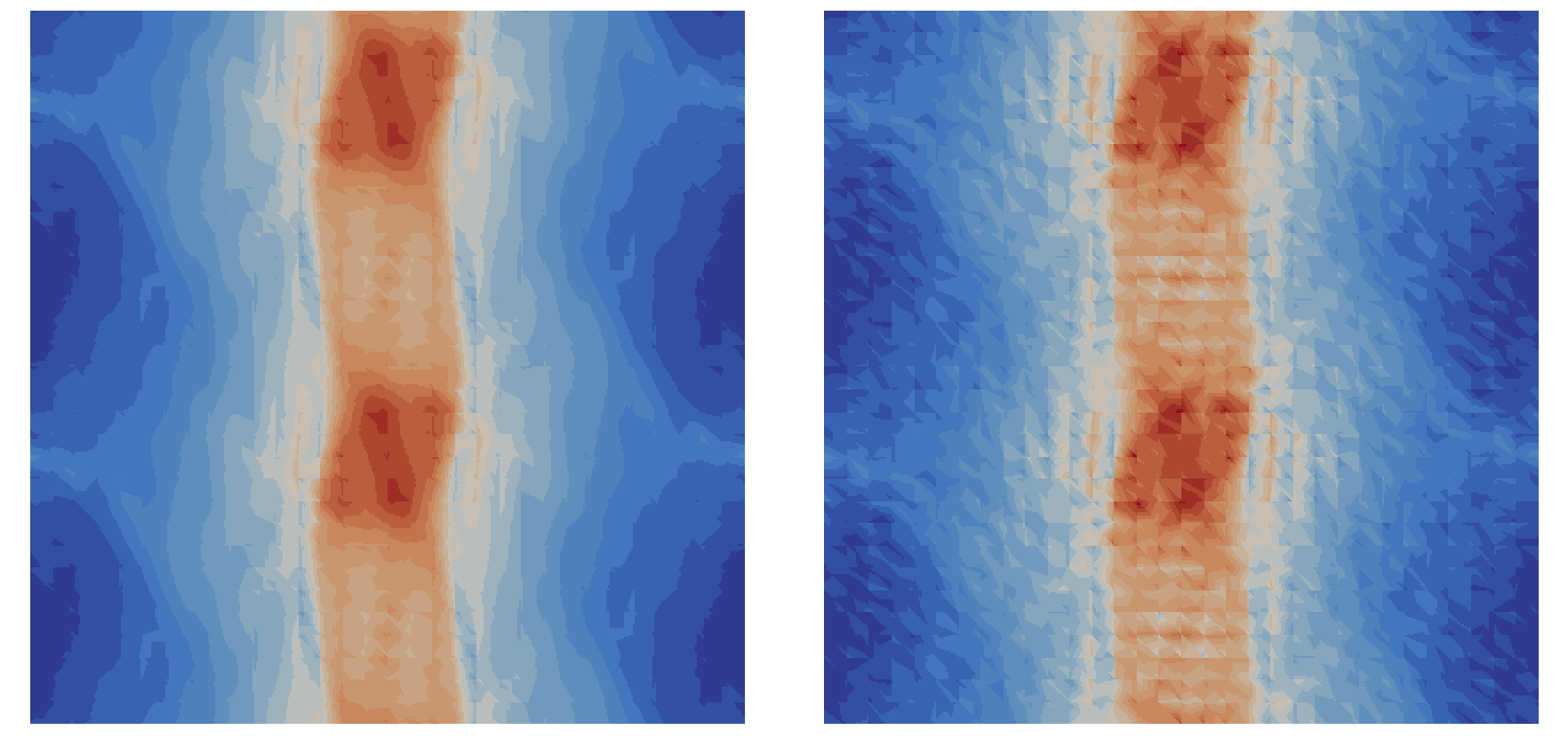}
\includegraphics[width=0.49\textwidth]{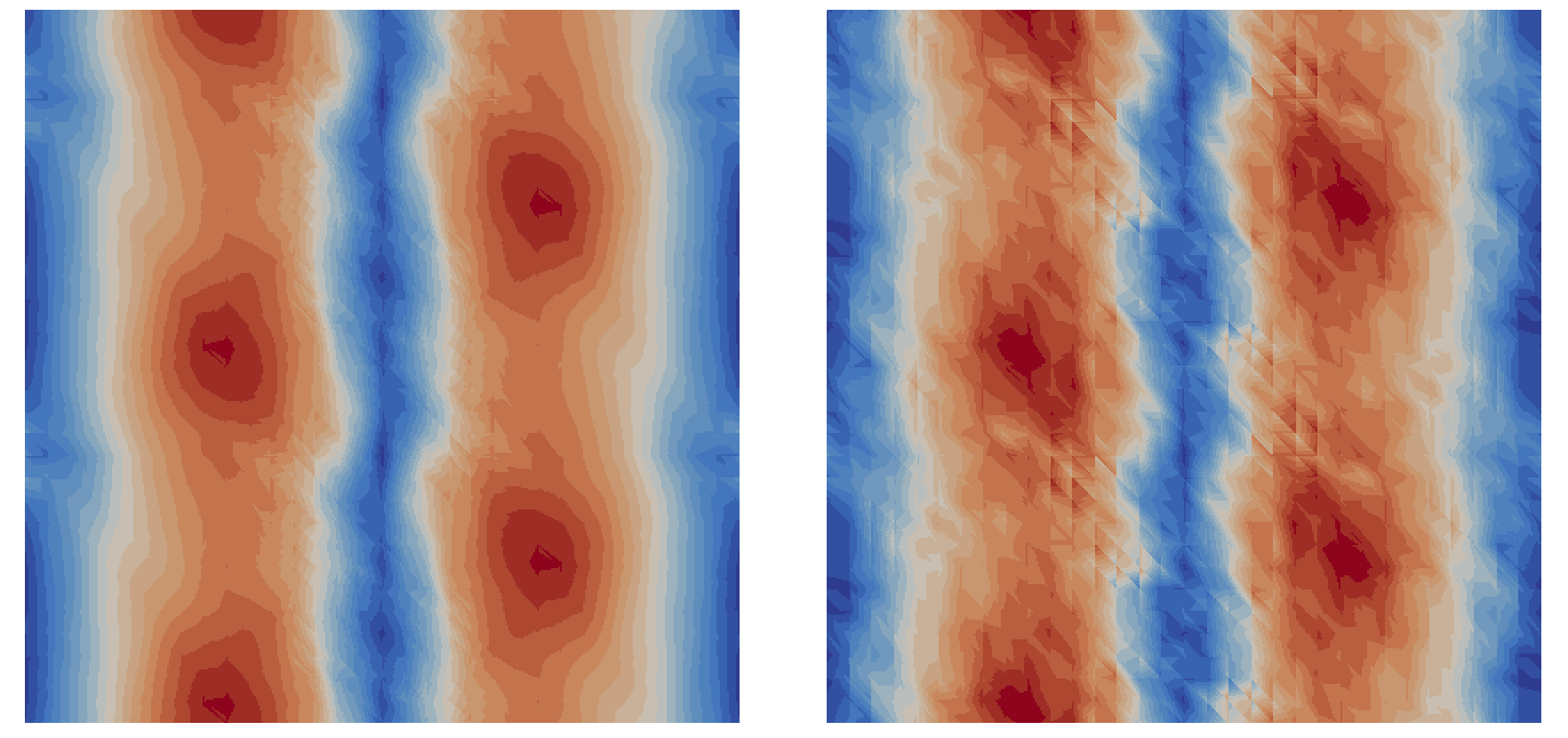}
\caption{Field development after 1000 time steps for periodic unit square test case with energy conserving setup including upwinding in $\mathbf{u}$. Left two images: depth fields, including upwinding and not including upwinding in $D$, respectively. Right two images: velocity fields, including upwinding and not including upwinding in $D$, respectively. Spatial resolution $32\times32$, $\Delta t =0.001$, with 4 Picard iterations for each time step. Depth field values 0.75 to 1.5 with contours every 0.05, velocity field magnitude values 0 to 1 with contours every 0.05.} \label{PUSM}
\end{center}
\end{figure}

\begin{figure}[ht]
\begin{center}
\includegraphics[width=1\textwidth]{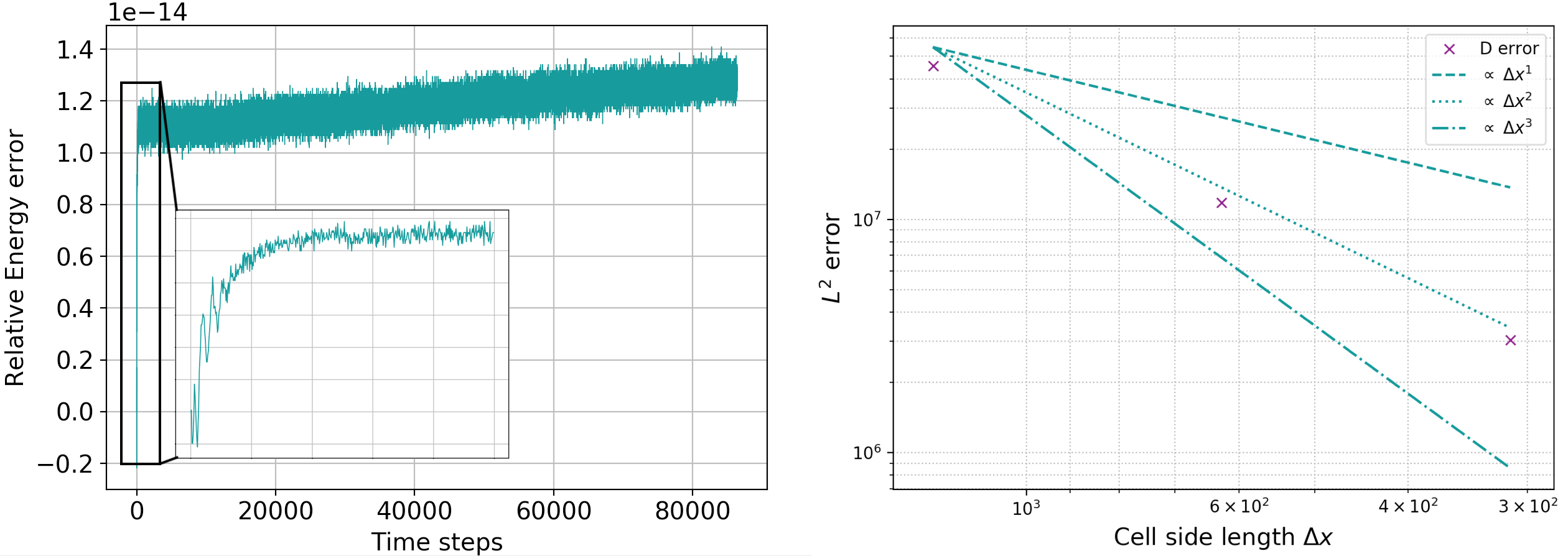}
\caption{Left: Relative energy error development for Williamson 2 test case, using fully energy conserving upwinded discretisation \eqref{discrete_u_eqn} - \eqref{discrete_D_eqn}, mesh refinement level 5, $\Delta t =50$s, with 4 Picard iterations for each time step (window: first 500 time steps). Right: L2 depth field error averaged over the last 1000 time steps for refinement levels 3 to 5.} \label{full_energy_conservation}
\end{center}
\end{figure}

\begin{figure}[ht]
\begin{center}
\includegraphics[width=1.0\textwidth]{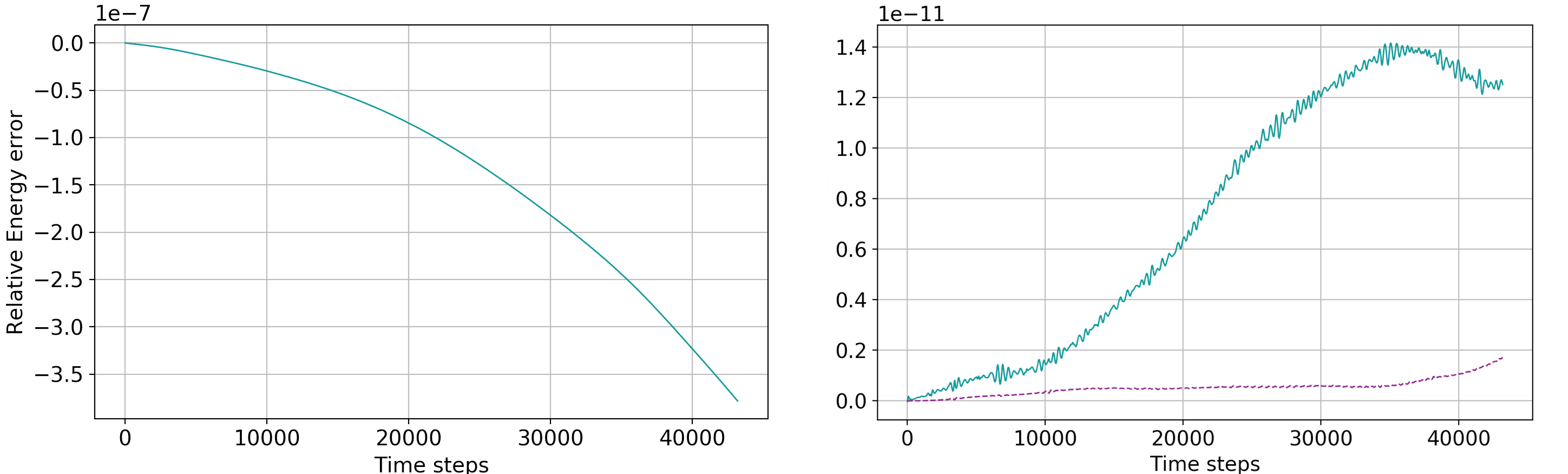}
\caption{Relative energy error developments for Williamson 5 test case. Left: non-energy conserving setup. Right: Energy conserving setup with (cyan) and without (dashed purple) upwinding in $D$. Mesh refinement level 5, $\Delta t =50$s, with 8 Picard iterations per time step.} \label{W5_energies}
\end{center}
\end{figure}

\begin{figure}[ht]
\begin{center}
\includegraphics[width=1.0\textwidth]{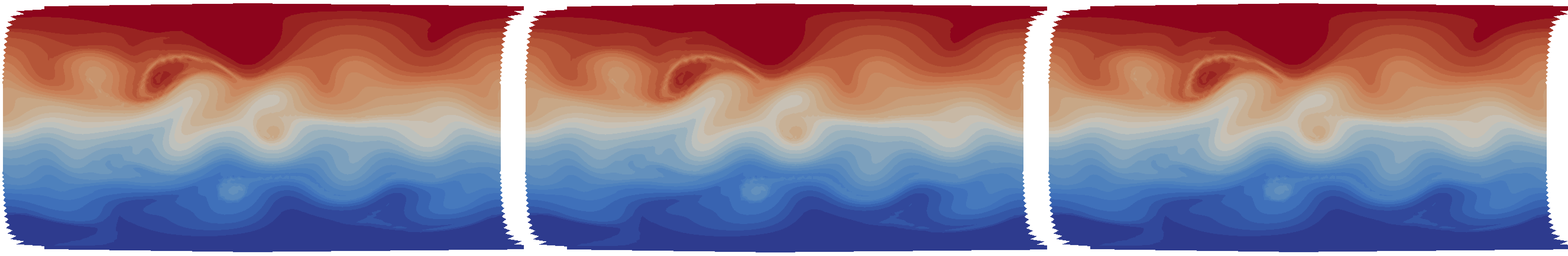}
\caption{Potential vorticity fields after 25 days for Williamson 5 test case. Left to right: Non-energy conserving setup, energy conserving setup with upwinding in $D$, energy conserving setup without upwinding in $D$. Mesh refinement level 5, $\Delta t =50$s, with 8 Picard iterations per time step. 30 contours, scale: $-3 \times 10^{-8}$ (blue) to $3 \times 10^{-8}$ (red).} \label{W5_fields}
\end{center}
\end{figure}

\begin{figure}[ht]
\begin{center}
\includegraphics[width=1.0\textwidth]{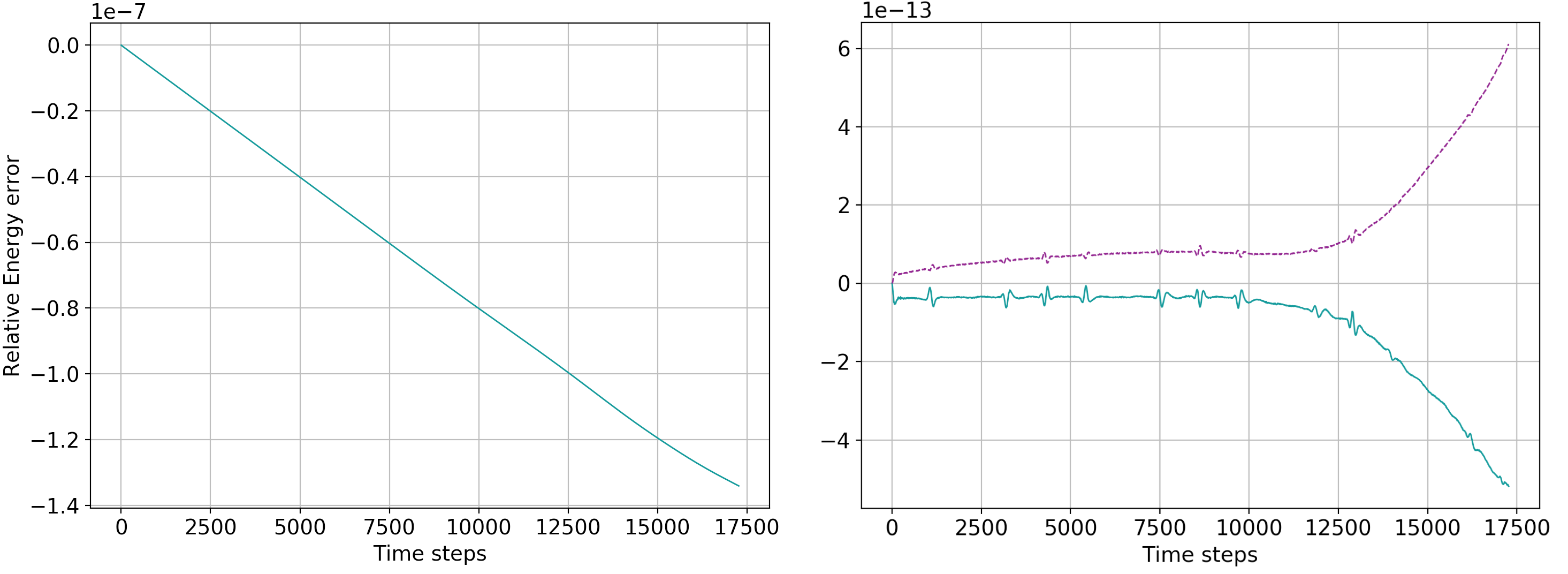}
\caption{Relative energy error developments for the Galewsky test case. Left: non-energy conserving setup. Right: Energy conserving setup with (cyan) and without (dashed purple) upwinding in $D$. Mesh refinement level 6, $\Delta t =30$s, with 8 Picard iterations for each time step.} \label{Galewsky_energies}
\end{center}
\end{figure}

\begin{figure}[ht]
\begin{center}
\includegraphics[width=1.0\textwidth]{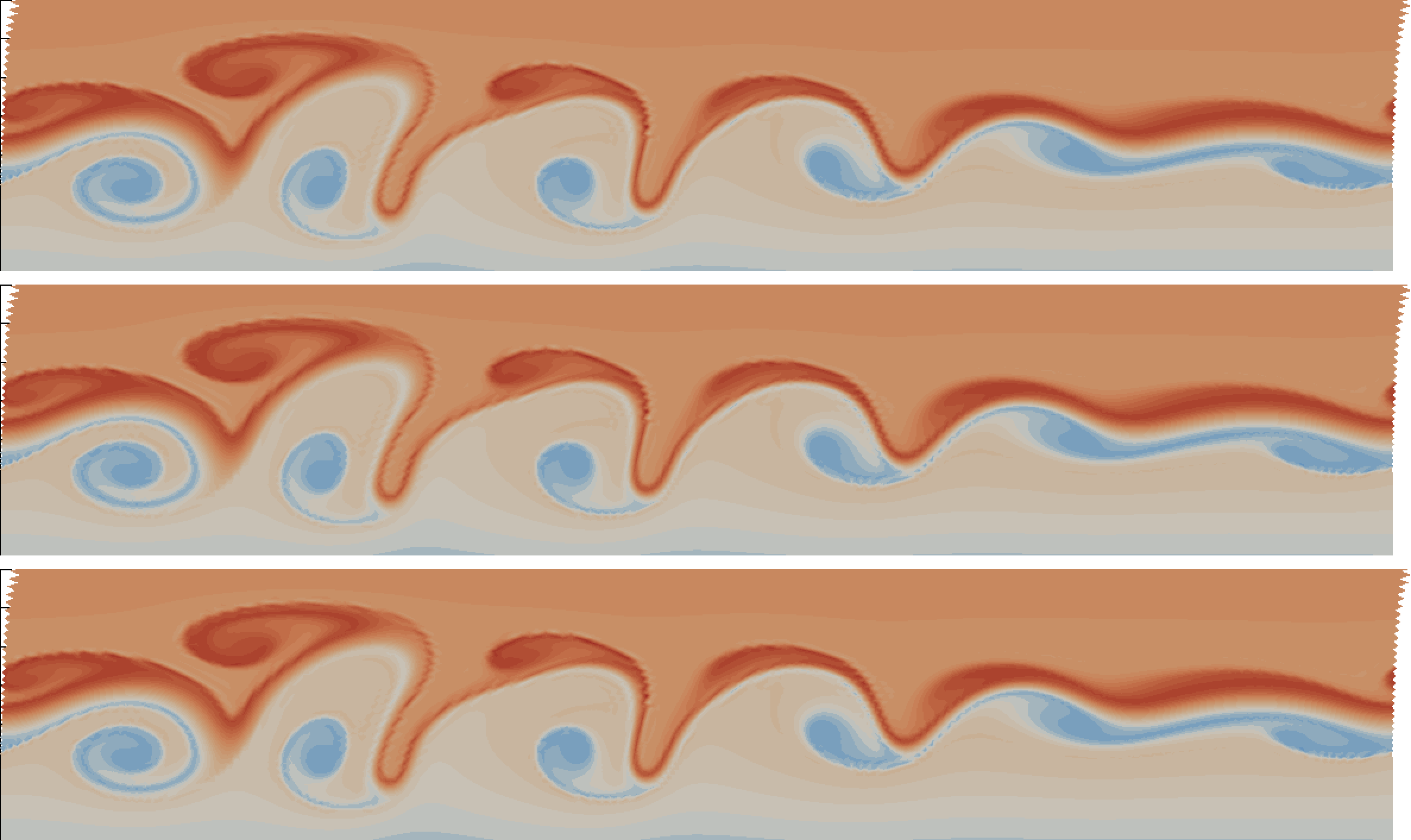}
\caption{Potential vorticity fields between 10 and 80 degrees latitude after 6 days for the Galewsky test case. Top to bottom: Non-energy conserving setup, energy conserving setup with upwinding in $D$, and without upwinding in $D$. Mesh refinement level 6, $\Delta t =30$s, with 8 Picard iterations for each time step. Contour lines every $1.25 \times 10^{-9}$, with negative values in blue, positive ones in red.} \label{Galewsky_fields}
\end{center}
\end{figure}

\begin{figure}[ht]
\begin{center}
\includegraphics[width=0.9\textwidth]{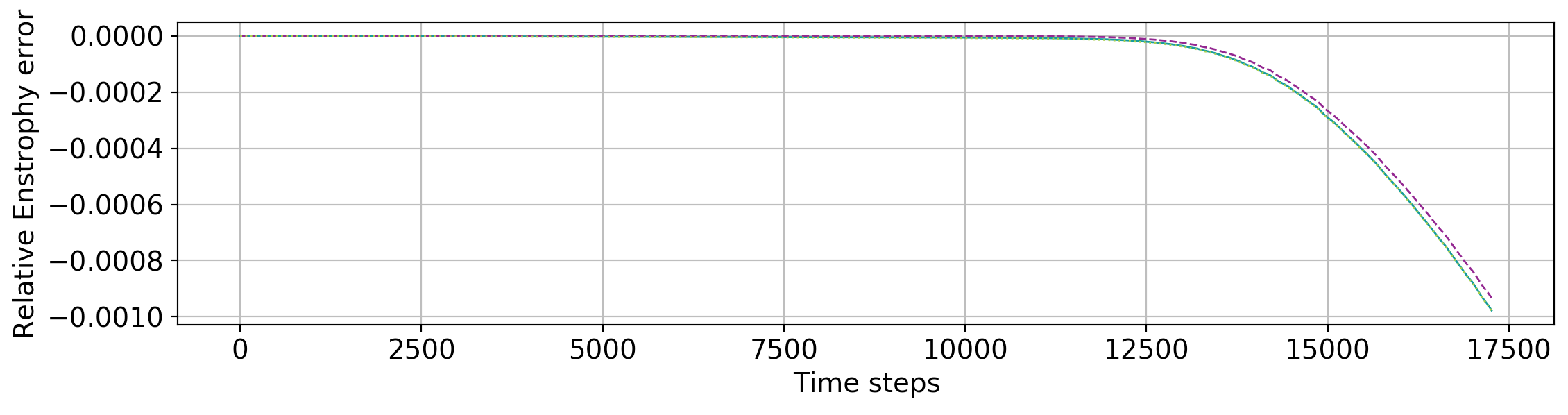}
\caption{Relative enstrophy error developments for the Galewsky test case, for energy conserving setup with (cyan) and without (dashed purple) upwinding in $D$ and non-energy conserving setup (dotted green).} \label{Galewsky_enstrophies}
\end{center}
\end{figure}

Having described the full discretisation as well as two other reference spatial discretisations, we proceed to our set of numerical tests. First, we consider a wave in a periodic unit square mesh as given in \cite{McRae_2014}. Being more of an artificial test case, it serves as a proof of concept for introducing upwinding in the velocity field $D$, and thus in extension for the density field $\rho$ in the context of the compressible Euler equations. Since we do not test for energy-conservation yet, we only consider the two energy-conserving versions, with and without upwinding in $D$. The initial conditions are given by
\begin{align}
\begin{split}
&\mathbf{u}_0 = (0, \; \sin(2 \pi x)), \label{unsteady_test_case}\\
&D_0 = 1 + \frac{1}{4 \pi} \frac{f}{g}\sin(4 \pi y),
\end{split}
\end{align}
with $(f, \; g) \coloneqq (5.0, \; 5.0)$. The domain is divided into $32^2$ squares, each of which in turn is divided into two triangles. The resulting fields for $\Delta t = 0.001$ and 1000 time steps, with 4 Picard iterations for each time step, are depicted in figure \ref{PUSM}.\\
\\
We find that the upwinding in $D$ not only significantly reduces small scale perturbations in the depth field, but also in the velocity field. Next, we consider more realistic spherical test cases. Since the depth field development in these cases is generally smoother than in the periodic unit square case, we anticipate the qualitative difference between the upwinded and non-upwinded versions to be small. Thus, the main purpose of these tests is to demonstrate an improved energy conservation as well as a qualitative field behaviour close to the projection-free non-energy conserving version.\\ \\
To validate the new upwinded scheme's energy conservation as well as consistency, we use the second of the standard Williamson spherical test cases given in \cite{williamson1992standard}, which describes a steady state scenario. Additionally, to compare the schemes with respect to their energy conservation properties as well as field development, we use the fifth test of the aforementioned test series, corresponding to flow past a mountain, as well as the Galewsky barotropic instability test case as described in \cite{galewsky2004initial}. In the Williamson 2 test case, the initial conditions are given by
\begin{align}
&\mathbf{u} = u_0(-y, x, 0)/a,\\
&D = h - (a\Omega u_0+ u_0^2/2)\frac{z^2}{ga^2},
\end{align}
for a sphere of radius $a = 6371220$m, with rotation rate $\Omega = 7.292 \times 10^{-5}$s$^{-1}$ (noting that $f = 2\Omega z /a$), and gravitational acceleration $g = 9.810616$ms$^{-2}$. The mean height and wind speed are given by $h = 5960$m and $u_0 = \frac{2\pi a}{12}$m/day. The simulation is run for 50 days, with a time step of $\Delta t = 50$s, and 4 Picard iterations for each time step. The mesh is given by an icosahedral triangulation, where refinement level 0 corresponds to 20 triangles. For every higher level, each triangle is refined to 4 triangles (so that each increase corresponds to halving the cell side length $\Delta x$). The resulting relative energy error development, as well as the L2 depth field error, averaged over the last 1000 time steps and for different refinement levels, are depicted in figure \ref{full_energy_conservation}.\\
\\
As expected, we find that energy is conserved up to machine precision throughout the simulation, with the initial increase likely due to the simplified variational derivative of the residual $\mathbf{R}$. Further, the L2 depth field error convergence as we refine the mesh also matches the expected second order rate (the same holding true for the velocity field).\\ \\
Next, we consider the fifth Williamson test, corresponding to unsteady flow over a mountain. The initial conditions are given by
\begin{align}
&\mathbf{u} = u_0(-y, x, 0)/a,\\
&D = h - (a\Omega u_0+ u_0^2/2)\frac{z^2}{ga^2} - b,\\
&b = b_0(1 - r/R),
\end{align}
where $b$ describes the mountain's surface, for $R=\pi/9$, mountain height $b_0 = 2000$m and $r$ such that $r = \text{min}(R, \sqrt{(\lambda - \lambda_c)^2 + (\theta - \theta_c)^2})$. $\lambda \in [-\pi, \pi]$ and $\theta \in [-\pi/2, \pi/2]$ denote longitude and latitude respectively, and the mountain's centre is chosen as $\lambda_c = -\pi/2$ and $\theta_c = \pi/6$. The mean height and wind speed are given by $h = 5960$m and $u_0 = 20$m/s. The simulation is run for 25 days, with a time step of $\Delta t = 50$s. Since this test involves an unsteady nonlinear field evolution, we increase the number of Picard iterations to 8 for each time step in order to keep the relative energy error contribution of the nonlinear solver small. The resulting relative energy errors and potential vorticity fields for the energy conserving setup with and without upwinding as well as the non-energy conserving setup are given in figures \ref{W5_energies} and \ref{W5_fields}. \\
\\
\\
Again, as expected the energy conserving formulations conserve energy up to a good degree, with a relative energy error of the order of $10^{-11}$, four orders of magnitude smaller than the relative energy error for the non-energy conserving setup. Additionally, we find a practically identical field development for all three setups, indicating that the new method for upwinding in $D$ in a Poisson bracket framework also behaves as expected in a spherical domain.
\begin{remark}
For 4 Picard iterations, the relative energy error for the energy conserving setups increases to the order of $10^{-9}$, while remaining at the order of $10^{-7}$ for the non-energy conserving space discretisation. In fact, at this stage the error originating from the nonlinear solver surpasses the decrease in error due to the energy-conserving time discretisation. We find that the midpoint rule time discretisation given in remark \ref{non_conservative_time_disc}, together with 4 Picard iterations, still yields relative energy errors of the order of $10^{-9}$ and $10^{-7}$ for the energy-conserving and non-energy conserving spatial discretisations, respectively. This demonstrates that the energy conserving spatial discretisations still lead to a significant improvement even for lower Picard iteration numbers and non-energy conserving time discretisations.
\end{remark}
Finally, we consider the Galewsky test case, simulating a barotropic instability. The initial conditions are given by a zonal flow $u$ confined within latitudes $\theta_0 = \pi/7$ and $\theta_1 = 5\pi/14$, and a background depth field $D$ in balance with $u$, perturbed by a localised bump $D_{p}$. The zonal flow is given by
\begin{equation}
u(\theta) = \frac{u_0}{e_n} \text{exp} \Big((\theta - \theta_0)(\theta - \theta_1)\Big)^{-1},
\end{equation}
for $u_0 = 80$m/s, and normalising constant $e_n = \text{exp}\big(-4/(\theta_1 - \theta_0)^2\big)$. To reach a steady state depth field, we integrate the continuity equation (in spherical form), leading to a field of form
\begin{equation}
gD(\theta) = g h_0 - \int^{\theta}_{-\frac{\pi}{2}} a u(\theta')\big(f + \frac{\text{tan}\theta'}{a} u(\theta')\big)d\theta',
\end{equation}
where $h_0$ is chosen such that the global mean height is equal to 10km. Finally, the perturbation is given by
\begin{equation}
D_p = h_p \text{cos}(\theta) \text{exp} \Big(-(\lambda/\alpha)^2 - ((\theta_2 - \theta)/\beta)^2\Big),
\end{equation}
for $h_p = 120$m , $\alpha = 1/3$, $\beta = 1/15$ and $\theta_2 = \pi/4$ (such that the perturbation is located directly in the zonal flow). The simulation is run for 6 days at mesh refinement level 6, with a time step of $\Delta t = 30$s, and 8 Picard iterations for each time step. The resulting relative energy errors as well as vorticity fields for the three setups are depicted in figures \ref{Galewsky_energies} and \ref{Galewsky_fields}.\\
\\
\\
We find a similar behaviour to the Williamson 5 test case. With respect to the relative energy error,  the two energy conserving setups outperform the non-energy conserving space discretisation by 6 orders of magnitude. Again, the field development is virtually the same for all setups.
\begin{remark}
Next to energy, the non-upwinded space discretisation \eqref{var_original_eqns} also conserves enstrophy. However, a controlled dissipation of enstrophy may be desirable, since the enstrophy cascades to small scales, eventually accumulating at the grid scale. This effect was countered in \cite{BAUER2018171} by including an SUPG scheme for the potential vorticity, which implies enstrophy dissipation for a sufficiently large SUPG parameter $\tau$. In our case, the upwinding in the velocity field $\mathbf{u}$ also dissipates enstrophy, with little difference in the dissipation rate whether or not upwinding in $D$ is also included in the discretisation (see image \ref{Galewsky_enstrophies}). More details on the dissipation of enstrophy depending on the choice of upwinding can be found in \cite{natale2017scale}.
\end{remark}
\section{Conclusion} \label{Conclusion}
In this paper, we introduced an energy conserving space discretisation for the rotating shallow water equations that includes upwinding in the depth and velocity fields. It is formulated using the compatible finite element method, and relies on a Hamiltonian framework with Poisson brackets to achieve energy conservation. The bracket is based on one without upwinding, which is described in \cite{McRae_2014}, and uses upwinding for the velocity field as formulated for the incompressible Euler equations in \cite{natale2016variational}. Upwinding for the depth field in this context has been newly introduced here, and relies on the introduction of an additional operator to recover the velocity field from the Hamiltonian variation corresponding to the momentum flux. In our numerical tests, we confirmed the scheme's energy conservation property, with a relative energy error close to machine precision when coupled with an energy conserving time discretisation. In the spherical test cases, we showed that the fully upwinded energy conserving scheme behaves as expected for spherical domains, leading to a field development virtually identical to that of a corresponding non-energy conserving upwinded reference scheme, despite the additional projections that are required to achieve energy conservation. Further, as demonstrated in the unit square test case, in the presence of large depth field gradients, the newly introduced upwinding in the depth field improves the field development compared to when upwinding is only applied for velocity, reducing small scale oscillations both in the depth and velocity fields.
\\ \\
The introduction of upwinding in the depth field was motivated by the development of a fully upwinded energy conserving space discretisation for the compressible Euler equations. In ongoing work, we aim to extend the Hamiltonian formulation presented here to the latter set of equations, incorporating upwinding for density as presented for the depth field in this paper, and additionally an SUPG formulation for the potential temperature field. While typical spherical shallow water scenarios feature a relatively small gradient in the depth field, the Euler equations often exhibit much stronger such gradients, therefore benefiting significantly from upwinding in the density field.
\newpage
\bibliography{SWE_var_upwind}
\end{document}